\documentclass[12pt]{article}
\usepackage{amsfonts,amssymb,amsmath,url}
\usepackage{xcolor}

\usepackage{amsthm}

\usepackage{amssymb}
\usepackage{graphics}
\usepackage{tikz}
\usetikzlibrary{shapes,backgrounds,calc}
\usepackage{latexsym}
\usepackage{multicol}
\usepackage{verbatim,enumerate}
\usepackage{accents}
\usepackage{cite}

\usepackage[colorlinks=true, linkcolor=blue, citecolor=blue, urlcolor=blue]{hyperref}
\usepackage{hyperref}
\usepackage{amsmath, amscd,url}

\usepackage{setspace}
\usepackage{pstricks}

\advance\textwidth by 1.3in \advance\oddsidemargin by -.6in \advance\evensidemargin by -.6in
\newtheorem{theorem}{Theorem}[section]
\newtheorem{prop}[theorem]{Proposition}
\newtheorem{lemma}[theorem]{Lemma}

\newtheorem{preprob}{Problem}
\newenvironment{prob}{\preprob\rm}{\endpreprob}

\newcommand{\grtype}{\mathop{\mbox{$\mathrm{A}$}}}
\newcommand{\partype}{\mathop{\mbox{$\mathrm{B}$}}}

\begin{document}

\title{Super graphs on groups, II}
\author{G. Arunkumar\\
{\small Department of Mathematics,  Indian Institute of Technology
Madras, Chennai - 600036, India}\\
Peter J. Cameron\\
{\small School of Mathematics and Statistics, University of St Andrews, Fife,
UK}\\
Rajat Kanti Nath\\
{\small Department of Mathematical Sciences, Tezpur University, Sonitpur, Assam 784028, India.}
}
\date{}
\maketitle

\begin{abstract}
In an earlier paper, the authors considered three types of graphs,
and three equivalence relations, defined on a group, viz.\ the power graph,
enhanced power graph, and commuting graph, and the relations of equality,
conjugacy, and same order; for each choice of a graph type $\grtype$ and
an equivalence relation $\partype$, there is a graph, the \emph{$\partype$ 
super$\grtype$ graph} defined on $G$. The resulting nine graphs (of which
eight were shown to be in general distinct) form a two-dimensional hierarchy.
In the present paper, we consider these graphs further. We  prove universality properties for the conjugacy
supergraphs of various types, adding the nilpotent,  solvable and enhanced power graphs to the commuting graphs considered in the rest of the paper, and also examine their relation to the
invariable generating graph of the group. We also show that supergraphs can
be expressed as graph compositions, in the sense of Schwenk, and use this
representation to calculate their Wiener index.
We illustrate these by computing Wiener index 
of
equality  supercommuting  and conjugacy supercommuting graphs for   dihedral and quaternion groups. 
\end{abstract}

\section{Introduction}
We begin with the definition of the types of graphs that we will be considering. We have
three basic types of graphs on a group $G$, defined by the condition for
$g$ and $h$ to be adjacent, as follows:
\begin{itemize}
\item the \emph{power graph}: one of $g$ and $h$ is a power of the other;
\item the \emph{enhanced power graph}: both $g$ and $h$ are powers of an
element $k$;
\item the \emph{commuting graph}: $gh=hg$.
\end{itemize}
Each  of the above class of graphs has a substantial literature. We refer to the recent survey
\cite{kscc} of power graphs, and the paper \cite{gong} for an account for the
hierarchy of graphs. In Section \ref{s:univ}, we also consider two
further graphs, with adjacency as follows:
\begin{itemize}
\item the \emph{nilpotent graph}: $\langle g,h\rangle$ is a nilpotent group;
\item the \emph{solvable graph}: $\langle g,h\rangle$ is a solvable group.
\end{itemize}

Also, we consider three equivalence relations, and corresponding partitions,
as follows:
\begin{itemize}
\item equality;
\item conjugacy;
\item same order.
\end{itemize}

For each graph type $\grtype$ and partition $\partype$, we define the
	\emph{$\partype$ super$\grtype$ graph} on a group $G$ by the rule that
	$g$ and $h$ are joined if and only if there exist $g'$ in $[g]_B$ and $h'$ in $[h]_B$ such that $g'$ is joined to 
	$h'$ in $\grtype(G)$. {\it By convention, we assume that any two vertices in the
same partition of $\partype$ are joined}.

Our definition would give the power graph the name ``equality superpower
graph'', but we will simply say ``power graph'', with similar convention for
the other basic graph types with the relation of equality.

This defines nine graphs. It is shown in \cite{acns} that they form a
two-dimensional hierarchy, in the sense
\begin{itemize}
\item for a given partition type $\partype$, the graphs are ordered by the
spanning subgraph relation in the order power graph, enhanced power graph,
and commuting graph;
\item for a given graph type $\grtype$, the graphs are ordered by the spanning
subgraph relation in the order equality, conjugacy, same order.
\end{itemize}
It is further shown that, for any group $G$, the order superenhanced power graph
and order supercommuting graph coincide, but no other pairs coincide in
general.

In this paper, we prove universality properties for the conjugacy
supergraphs of various types consisting of enhanced power graph, commuting graph, nilpotent graph and  solvable graph. A relation between conjugacy super$\grtype$ and invariable generating graph is also established.  
We will  express the supergraphs as graph compositions, in the sense of Schwenk \cite{sch},  and
use this expression to examine the Wiener index \cite{w}
in the case of some special groups. See \cite{survey, survey1} and their big list of references for the literature on the Weiner index of graphs. The spectral aspects of these graphs are discussed in \cite{ACN-partIII}.


\section{Universality of supergraphs}
\label{s:univ}

In what follows, we use the graph-theoretic notions from \cite{dbw}. A class of finite graphs is called \emph{universal} if every finite graph can
be embedded as  an  induced subgraph in some graph in the class.

Universal graphs were studied in  \cite{egg21, pr99} . In \cite{cg83}, the search for a graph $G$ on $n$ vertices with the minimum number of edges such that every tree on $n$ vertices is isomorphic to a spanning tree of G is carried out. 
In 
\cite{gong}, it is shown that various graph types, including the enhanced
power graph, commuting graph, nilpotent graph and solvable graph, on
finite groups form universal classes. Similarly, in the recent paper \cite{ackt23},  the following results are proved. 

\begin{enumerate}
	\item The zero-divisor graphs of Boolean rings are universal. That is, every finite graph is an induced subgraph of the zero-divisor graph of a Boolean ring. 
	
	\item  The class of zero-divisor graphs of the rings of integers modulo n is universal. Further, we can even restrict the $n$ to be squarefree. 
	
	\item For every finite graph $G$, there is a finite commutative local ring R with unity such that $G$ is an induced subgraph of the zero-divisor graph of R. 
\end{enumerate}

Here we examine universality for
conjugacy supergraphs. 

Recall the \emph{condensed} versions of the conjugacy
supergraphs, where the vertices are the conjugacy classes of elements, two
classes $x^G$ and $y^G$ being joined if and only if there exist $x'\in x^G$
and $y'\in y^G$ such that $x'\sim y'$ ($x'$ is adjacent to $y'$) in the original (uncompressed) graph.
The compressed conjugacy supergraphs for the commuting, nilpotent and
solvable graphs have been studied under the names
commuting/nilpotent/solvable conjugacy class graph, and denoted $\Gamma_{ccc}$,
$\Gamma_{ncc}$ and $\Gamma_{scc}$ respectively, see \cite{hlm,ME2017,bcns}.

Power graphs of finite groups do not form a universal graph, since any power
graph is the comparability graph of a partial order. The same applies to
the conjugacy superpower graphs. But the following holds:

\begin{prop}
A finite graph $\mathcal{G}$ is embeddable in the conjugacy superpower graph
(or its compressed version) if and only if it is the comparability graph of
a finite partial order.
\end{prop}

\begin{proof}
The corresponding result for the power graph is \cite[Theorem 5.4]{gong},
where the group used for the embedding is abelian (indeed, cyclic of squarefree
order); every conjugacy class is a singleton, so the power graph and the
compressed and uncompressed conjugacy superpower graphs all coincide. 
\end{proof}

\begin{theorem}
Let $\mathcal{A}$ be the set of graph types consisting of enhanced power graph,
commuting graph, nilpotent graph, and solvable graph. Then, for any
$\grtype$ in $\mathcal{A}$, the classes of uncompressed and compressed
conjugacy super$\grtype$ graphs are universal.
\end{theorem}

\begin{proof}
We begin by noting that the compressed graphs are induced subgraphs of the
uncompressed graphs (on a set of conjugacy class representatives), and so it
suffices to prove the theorem for the compressed graphs.

We will first deal with the commuting, nilpotent and solvable graphs, all
of which are handled by the same argument: given a graph $\Gamma$, we construct
an embedding of its vertices into a group $G$ such that, if two vertices $x^G$
and $y^G$ are adjacent, then there are elements $x'$ and $y'$ in these classes
such that $x'$ and $y'$ commute, whereas if they are non-adjacent, then
$\langle x',y'\rangle$ is non-solvable for any choice of $x'$ and $y'$ in these
classes. At the end, we indicate how to modify the argument for the enhanced
power graph.

So let $\Gamma(G)$ denote the condensed conjugacy super$\grtype$ graph of
$G$, where $\grtype$ denotes the commuting, nilpotent or solvable graph.

\paragraph{Step 1} The \emph{strong product} $\mathcal{G}\boxtimes\mathcal{H}$
of two graphs $\mathcal{G}$ and $\mathcal{H}$ is the graph with vertex set
$V(\mathcal{G})\times V(\mathcal{H})$ (where $V(\mathcal{G})$ is the vertex set of $\mathcal{G}$), and distinct vertices $(u,u{'})$ and $(v,v{'})$ are adjacent if and only if: \\
$u = v$ and $u{'}$ is adjacent to $v^{'}$, or \\
$u{'} = v{'}$ and $u$ is adjacent to $v$, or \\
$u$ is adjacent to $v$ and $u{'}$ is adjacent to $v{'}$.

The symbol $\boxtimes$ is intended to suggest the
strong product of two copies of $K_2$.

Now we claim that
\[\Gamma(G\times H)=\Gamma(G)\boxtimes\Gamma(H).\]
To see this, take $g_1,g_2\in G$ and $h_1,h_2\in H$. Then 
$\langle(g_1,h_1),(g_2,h_2)\rangle$ is a subgroup of
$\langle g_1,g_2\rangle\times\langle h_1,h_2\rangle$, and projects onto each
factor of this product (that is, it is a subdirect product); so it is
abelian/nilpotent/solvable if and only if the two factors are.

\paragraph{Step 2} If $\mathcal{G}$ and $\mathcal{H}$ are graphs on the same
vertex set $V$, then we define their intersection $\mathcal{G}\cap\mathcal{H}$
to have the same vertex set as each of the graphs, its edge set being the
intersection of the edge sets.

Now observe that $\mathcal{G}\cap\mathcal{H}$ is an induced subgraph of
$\mathcal{G}\boxtimes\mathcal{H}$, corresponding to the vertices $(x,x)$ on the
diagonal of the cartesian product.

\paragraph{Step 3} Given $n$, the complete graph $K_n$, and the graph
$K_n-e$ obtained by removing one edge, has the property that its vertex set
is embeddable in some group $G$ in such a way that adjacent vertices commute
while non-adjacent vertices generate a non-solvable group. We show this for
$K_n-e$; the proof for $K_n$ is similar but easier. Take vertex set
$\{1,\ldots,n\}$, and let $e=\{n-1,n\}$. We may clearly assume that $n\ge3$.

Take $p_1,p_2,\ldots,p_n$ to be the
first $n$ primes in order. Let $N=p_{n-1}+p_n-1$, and let $G$ be the symmetric
group of degree $N$. Take $g_i$ to be a $p_i$-cycle in $G$. Then, if
$\{i,j\}\ne\{n-1,n\}$, then $p_i+p_j\le N$, and we can find disjoint (and
hence commmuting) cycles of lengths $p_i$ and $p_j$; but $p_{n-1}+p_n>N$,
so any two cycles of lengths $p_{n-1}$ and $p_n$ intersect, and generate an
alternating group.

\paragraph{Step 4} Proof of the theorem. Given a graph $\mathcal{G}$ on $n$
vertices, we can assume that $\mathcal{G}$ is not complete; write it as the
intersection of the graphs $K_n-e$ where $e$ runs over all \emph{nonedges} of
$\mathcal{G}$. By Step~3, all these graphs are embeddable in compressed
conjugacy supergraphs of any of the three types, in groups $G_e$; then
$\mathcal{G}$ is embedded in the
corresponding graph for the direct product of the groups $G_e$.

\paragraph{The enhanced power graph} requires a different argument; for the
direct product of two cyclic groups is cyclic if and only if the groups have
coprime orders; but given a family of groups of pairwise coprime orders, all
but one have odd order, and so are solvable by the Odd-Order Theorem. So
we proceed a little differently. For each edge of the graph $\mathcal{G}$, we
choose a set of $n$ primes, so that these sets are pairwise disjoint. If the
primes for a particular edge $e$ are $p_1^e,\ldots,p_n^e$ in ascending order,
we take the symmetric group of degree $N_e=p_n^e+p_{n-1}^e-1$, and represent
vertices  of the graph by the conjugacy classes of cycles of length
$p_1^e,\ldots,p_n^e$ as before. This gives an embedding of $K_n-e_0$, where
$e_0=\{n-1,n\}$. Now, in the direct products of these groups, the cycles
representing a vertex $i$ in the different factors have distinct prime lengths;
so, if these cycles commute, they generate a cyclic
group whose order is the product of these primes. 
\end{proof}

\section{The invariable generating graph}

In this section, we discuss a graph which has begun to attract attention very
recently, and show how it is related to conjugacy supergraphs.

Let $G$ be a finite group. A set $\{x_1,\ldots,x_n\}$ of elements of $G$ is
said to \emph{invariably generate} $G$ if, for any $g_1,\ldots,g_n\in G$,
the set $\{x_1^{g_1},\ldots,x_n^{g_n}\}$ generates $G$; in other words,
the original elements can be replaced by arbitrary conjugates and we still
have a generating set.

It is known that any finite simple group can be invariably generated by two
elements~\cite{gm,kls}.

The \emph{generating graph} of $G$ is the graph with vertex set $G$, in which
two elements are joined if they generate $G$. By analogy, we define the
\emph{invariable generating graph} of $G$ to be the graph with vertex set $G$,
in which two elements are joined by an edge if and only if they invariably
generate $G$.

For the generating graph, the following holds
\cite[Propositions~2.6,~11.1]{gong}:

\begin{prop}
	\begin{enumerate}
		\item
		If $G$ is non-abelian, then the generating graph is contained in the complement
		of the commuting graph, with equality if and only if $G$ is a minimal
		non-abelian group.
		\item
		If $G$ is non-nilpotent, then the generating graph is contained in the
		complement of the nilpotent graph, with equality if and only if $G$ is a
		minimal non-nilpotent group.
		\item
		If $G$ is non-solvable, then the generating graph is contained in the
		complement of the solvable graph, with equality if and only if $G$ is a
		minimal non-solvable group.
	\end{enumerate}
	\label{p:nonsuper}
\end{prop}

The minimal non-abelian groups were classified by Miller and Moreno~\cite{mm},
and the minimal non-nilpotent groups by Schmidt~\cite{schmidt} (the latter are
often called \emph{Schmidt groups}). If $G$ is minimal non-solvable, then $G$
modulo its solvable radical is a minimal simple group, and occurs in the
classification of N-groups by Thompson~\cite{thompson}.

The invariable generating graph bears a similar relation to supergraphs.

\begin{prop}
	Let $\grtype$ denote one of the properties ``abelian'', ``nilpotent'', ``solvable''.
	If $G$ is non-$\grtype$, then the invariable generating graph is contained in
	the complement of the conjugacy super$\grtype$ graph.
	\label{p:inv-super}
\end{prop}

\begin{proof}
	Suppose that $G$ is non-abelian. If $x$ and $y$ are adjacent in the invariable
	generating graph, then for any conjugates $x'$ and $y'$ of $x$ and $y$
	respectively, $x'$ and $y'$ generate $G$, and so cannot commute; this $x$ and
	$y$ are non-adjacent in the commuting graph.
	
	The proofs in the other two cases are similar.
\end{proof}

However, it is not known which groups satisfy equality here. If the invariable
generating graph is the complement of the commuting graph, then $G$ has
the property that, if $x$ and $y$ do not generate $G$, then there are 
conjugates $x'$ and $y'$ of $x$ and $y$ respectively which commute.

\begin{prob}
	Find all groups satisfying equality in Proposition~\ref{p:inv-super}.
\end{prob}

\section{Graph composition and quotient graph}

In \cite{sch}, Schwenk introduced the following notion of generalized composition of graphs and  denoted it by $\mathcal{H}[\Gamma_1,\Gamma_2,\dots,\Gamma_k]$. 
(This operation has been studied under other names such as generalized
lexicographic product and joined union in \cite{accv,gerbaud,mwg,ww}.)


\paragraph{Definition} Let $\mathcal{H}$ be a graph with vertex set $V(\mathcal{H}) = \{1,\dots,k\}$ and let $\Gamma_1,\dots,\Gamma_k$ be a collection of graphs with the respective vertex sets $V(\Gamma_i) = \{v_i^{1},\dots,v_i^{n_i}\}$ for $(1 \le i \le k)$. Then their generalized composition $\mathcal{G} = \mathcal{H}[\Gamma_1,\dots,\Gamma_k]$ has vertex set $V(\Gamma_1) \sqcup \cdots \sqcup V(\Gamma_k)$ and two vertices $v_i^p$ and $v_j^q$ of $\mathcal{G}$ are adjacent if one the following conditions is satisfied:
\begin{enumerate}
	\item $i=j$, and $v_i^p$ and $v_j^q$ are adjacent vertices in $\Gamma_i$.
	\item $i \ne j$ and $i$ and $j$ are adjacent in $\mathcal{H}$. 
\end{enumerate}

The graph $\mathcal{H}$ is said to be the base graph of $\mathcal{G}$ and the
graphs $\Gamma_i$ are called the factors of $\mathcal{G}$.

We observe that the join $\Gamma_1 \vee \Gamma_2$ of two graphs $\Gamma_1$ and $\Gamma_2$  is same as the generalized composition $K_2[\Gamma_1,\Gamma_2]$ where $K_2$ is the complete graph on two vertices.

\begin{prop}\label{p:join}
	Let $\grtype$ be a graph and $\partype$ a partition on a group $G$. Then the
	$\partype$ super$\grtype$ graph on $G$ is isomorphic to a generalized
	composition $\Delta[K_{n_1},\ldots,K_{n_r}]$ for some graph $\Delta$ on $r$
	vertices, where $n_1,\ldots,n_r$ are the sizes of the equivalence classes
	of $\partype$.
\end{prop}

The graph $\Delta$ is the induced subgraph on the set of equivalence class
representatives of $\partype$. We note that, if we had not made the convention
about the equivalence classes inducing complete subgraphs, we would take the
induced subgraphs on the parts of the partition.

The graph $\Delta$ will be called the $\partype$ quotient-super$\grtype$ graph
on $G$. It may be much smaller than $\grtype(G)$. For example, if $\partype$
is the ``same order'' relation, then the vertices can be taken to be the
orders of elements of $G$; if $\grtype$ is the power graph, then two vertices
are joined if and only if one divides the other, as explained earlier.

The expression of one of our supergraphs in terms of generalized composition,
and knowledge of the quotient graph, are useful in some contexts. Later in
the paper we apply these to the computation of the Wiener index and the
spectrum. However, when comparing different graphs in the hierarchy, it is
necessary to have a common vertex set, so we would not usually use these
tools.

\section{Realizing supergraphs}
In this section we realize  equality supercommuting  and conjugacy supercommuting graphs 
of  dihedral grups and  generalized quaternion groups. Following \cite{acns}, we use the notations $\mathrm{ESCom}(G)$ and $\mathrm{CSCom}(G)$ to denote the equality supercommuting and conjugacy supercommuting graphs of a group $G$.  Note that  equal supercommuting graph of a group coincides with the commuting graph of that group with the identification $\{g\} \leftrightarrow g$ between the vertices. 
\begin{theorem}\label{ESCom-D2n}
	Let   $D_{2n} = \langle a, b : a^n =b^2= e, bab^{-1} = a^{-1}\rangle$ be the dihedral group. Then
	\[
	\mathrm{ESCom}(D_{2n}) \cong \begin{cases}
		K_1 \vee (\underbrace{K_1\sqcup \cdots \sqcup K_1}_{n\text{-times}} \sqcup K_{n - 1}), &\text{ if  $n$ is odd}\\
		K_2 \vee (\underbrace{K_2\sqcup \cdots \sqcup K_2}_{\frac{n}{2}\text{-times}} \sqcup K_{n - 2}), &	\text{ if $n$ is even.}
	\end{cases}
	\]
\end{theorem}
\begin{proof}
If $n$ is odd then the  centralizers of elements of $D_{2n}$ are given by
\[
C_{D_{2n}}(e) = D_{2n}, C_{D_{2n}}(a^i) = \langle a \rangle \text{ for } 1 \leq i \leq n - 1
\] 
and 
\[
C_{D_{2n}}(a^ib) =\{e, a^ib\} \text{ for } 1 \leq i \leq n. 
\]
In $\mathrm{ESCom}(D_{2n})$, $e$ is the only dominant vertex  and it is  the only vertex adjacent to $a^ib$ for $1 \leq i \leq n$.  Let $\Gamma[X]$ be the subgraph of any graph $\Gamma$ induced by the subset $X$ of $V(\Gamma)$. Then $\mathrm{ESCom}(D_{2n})[S_1] \cong \underbrace{K_1\sqcup \cdots \sqcup K_1}_{n\text{-times}}$, where  $S_1 = \{b, ab, \dots, a^{n-1}b\}$. Also,
$\mathrm{ESCom}(D_{2n})[S_2]$ is complete if $S_2=\{e, a, a^2, \dots, a^{n - 1}\}$. It follows that 
\begin{align*}
	&\mathrm{ESCom}(D_{2n}) \\
	&~~~~~~= \mathrm{ESCom}(D_{2n})[\{e\}] \vee (\mathrm{ESCom}(D_{2n})[S_1] \sqcup \mathrm{ESCom}(D_{2n})[S_2 \setminus \{e\}])\\
	&~~~~~~\cong K_1 \vee (\underbrace{K_1\sqcup \cdots \sqcup K_1}_{n\text{-times}} \sqcup K_{n - 1}),
\end{align*}
since $\mathrm{ESCom}(D_{2n})[\{e\}] \cong K_1$ and $\mathrm{ESCom}(D_{2n})[S_2 \setminus \{e\}] \cong  K_{n - 1}$.	

If $n$ is even then  centralizers of elements of $D_{2n}$ are given by
\[
C_{D_{2n}}(e) = D_{2n} = C_{D_{2n}}(a^{\frac{n}{2}}), C_{D_{2n}}(a^i) = \langle a \rangle \text{ for } 1 \leq i \leq n - 1 \text{ and } i \ne \frac{n}{2}
\]
and
\[  
C_{D_{2n}}(a^ib) =\{e, a^{\frac{n}{2}},  a^ib, a^{i + \frac{n}{2}}b\} \text{ for } 1 \leq i \leq n.
\]
In $\mathrm{ESCom}(D_{2n})$, $e$ and $a^{\frac{n}{2}}$ are  the only dominant vertices. The only vertex, other than the dominant vertices,  adjacent to $a^ib$ is $a^{i + \frac{n}{2}}b$ for $1 \leq i \leq n$. Therefore,    
$\mathrm{ESCom}(D_{2n})[\{a^ib, a^{i + \frac{n}{2}}b\}] \cong K_2$ for $1 \leq i \leq n$ and $\mathrm{ESCom}(D_{2n})[S_1] \cong \underbrace{K_2\sqcup \cdots \sqcup K_2}_{\frac{n}{2}\text{-times}}$, where $S_1 = \{b, a^{\frac{n}{2}}b, ab, a^{\frac{n}{2} + 1}b, \dots, a^{\frac{n}{2}-1}b,$  $a^{n-1}b\}$. Also, 
$\mathrm{ESCom}(D_{2n})[S_2]$ is complete if $S_2=\{e, a^{\frac{n}{2}}, a, a^2, \dots, a^{\frac{n}{2} -1}$, $a^{\frac{n}{2}+1}, a^{n-1}\}$. It follows that 
\begin{align*}
	&\mathrm{ESCom}(D_{2n}) \\
	&~~= \mathrm{ESCom}(D_{2n})[\{e, a^{\frac{n}{2}}\}] \vee (\mathrm{ESCom}(D_{2n})[S_1] \sqcup \mathrm{ESCom}(D_{2n})[S_2 \setminus \{e, a^{\frac{n}{2}}\}])\\
	&~~\cong K_2 \vee (\underbrace{K_2\sqcup \cdots \sqcup K_2}_{\frac{n}{2}\text{-times}} \sqcup K_{n - 2}),
\end{align*}
since $\mathrm{ESCom}(D_{2n})[\{e, a^{\frac{n}{2}}\}] \cong K_2$ and $\mathrm{ESCom}(D_{2n})[S_2 \setminus \{e, a^{\frac{n}{2}}\}] \cong  K_{n - 2}$. 
\end{proof}

\begin{theorem}\label{ESCom-Q4n}
Let  $Q_{4n} = \langle a, b : a^{2n} = e, a^n =b^2, bab^{-1} = a^{-1}\rangle$ be the generalized quaternion group. Then
\[
\mathrm{ESCom}(Q_{4n}) \cong  K_2 \vee (\underbrace{K_2\sqcup \cdots \sqcup K_2}_{n\text{-times}} \sqcup K_{2n - 2}).
\]
\end{theorem}
\begin{proof}
The centralizers of elements of $Q_{4n}$ are given by
\[
C_{Q_{4n}}(e) = Q_{4n} = C_{Q_{4n}}(a^n), C_{Q_{4n}}(a^i) = \langle a \rangle \text{ for } 1 \leq i \leq 2n -1 \text{ and } i \ne n
\]
and
\[  
C_{Q_{4n}}(a^ib) =\{e, a^n,  a^ib, a^{i + n}b\} \text{ for } 1 \leq i \leq 2n. 
\]
In $\mathrm{ESCom}(Q_{4n})$, $e$ and $a^n$ are  the only dominant vertices. The only vertex, other than the dominant vertices,  adjacent to $a^ib$ is $a^{n + i}b$ for $1 \leq i \leq 2n$. Therefore,    
$\mathrm{ESCom}(Q_{4n})[\{a^ib, a^{n + i}b\}] \cong K_2$ for $1 \leq i \leq 2n$ and $\mathrm{ESCom}(Q_{4n})[S_1] \cong \underbrace{K_2\sqcup \cdots \sqcup K_2}_{n\text{-times}}$, where $S_1 = \{b, a^nb, ab, a^{n + 1}b, \dots, a^{n-1}b,$  $a^{2n-1}b\}$. Also, 
$\mathrm{ESCom}(Q_{4n})[S_2]$ is complete if $S_2=\{e, a^n, a, a^2, \dots, a^{n -1}$, $a^{n+1}, a^{2n-1}\}$. It follows that 
\begin{align*}
	&\mathrm{ESCom}(Q_{4n}) \\
	&~~= \mathrm{ESCom}(Q_{4n})[\{e, a^n\}] \vee (\mathrm{ESCom}(Q_{4n})[S_1] \sqcup \mathrm{ESCom}(Q_{4n})[S_2 \setminus \{e, a^n\}])\\
	&~~\cong K_2 \vee (\underbrace{K_2\sqcup \cdots \sqcup K_2}_{n\text{-times}} \sqcup K_{2n - 2}),
\end{align*}
since $\mathrm{ESCom}(Q_{4n})[\{e, a^n\}] \cong K_2$ and $\mathrm{ESCom}(Q_{4n})[S_2 \setminus \{e, a^n\}] \cong  K_{2n - 2}$.
\end{proof}

\begin{theorem}\label{CSCom-D2n}
	Let   $D_{2n} = \langle a, b : a^n =b^2= e, bab^{-1} = a^{-1}\rangle$ be the dihedral group. Then
	\[
	\mathrm{CSCom}(D_{2n}) \cong \begin{cases}
		K_1 \vee (K_1 \sqcup K_{\frac{n - 1}{2}})[K_1, \underbrace{K_2, \dots, K_2}_{(\frac{n-1}{2})\text{-times}}, K_n], \text{ if  $n$ is odd}\\
		K_2 \vee (K_1 \sqcup K_1 \sqcup K_{\frac{n}{2} - 1})[K_1, K_1, \underbrace{K_2, \dots, K_2}_{(\frac{n}{2}-1)\text{-times}}, K_{\frac{n}{2}}, K_{\frac{n}{2}}], 	\text{ if $n$ and $\frac{n}{2}$ are even}\\
		K_2 \vee (K_2 \sqcup K_{\frac{n}{2} - 1})[K_1, K_1, \underbrace{K_2, \dots, K_2}_{(\frac{n}{2}-1)\text{-times}}, K_{\frac{n}{2}}, K_{\frac{n}{2}}], 	\text{ if $n$ is even and $\frac{n}{2}$ is odd}
	\end{cases}
	\]
\end{theorem}
\begin{proof}
If $n$ is odd then the conjugacy classes and centralizers of elements of $D_{2n}$ are given by

$$e^{D_{2n}} = \{e\}, 
(a^j)^{D_{2n}} = \{a^j, a^{n - j}\} \text{ for } 1 \leq j \leq \frac{n - 1}{2},  b^{D_{2n}} = \{b, ab, \dots, a^{n - 1}b\}$$
and  
\begin{align*}
	C_{D_{2n}}(e) = D_{2n}, &C_{D_{2n}}(a^i) = \langle a \rangle \text{ for } 1 \leq i \leq n - 1\\  &C_{D_{2n}}(a^ib) =\{e, a^ib\} \text{ for } 1 \leq i \leq n, 
\end{align*}
\noindent respectively. It can be seen that 
$$\mathrm{CSCom}(D_{2n}) \cong \Delta[K_1, \underbrace{K_2, \dots, K_2}_{(\frac{n-1}{2})\text{-times}}, K_n],$$ 
where $V(\Delta) = \{e, a, a^2, \dots, a^{\frac{n - 1}{2}}, b\}$ which is identified with the set $\{1, 2, \dots,$  $\frac{n + 3}{2}\}$ preserving the order. Note that $e$ is the only dominant vertex of $\Delta$ and it is also the only vertex adjacent to $b$ (in $\Delta$). 
We have  $\Delta[S]$ is complete if $S=\{e, a, a^2, \dots, a^{\frac{n - 1}{2}}\}$.
Therefore 
$$\Delta = \Delta[\{e\}] \vee (\Delta[\{b\}] \sqcup \Delta[S \setminus \{e\}]) \cong K_1 \vee (K_1 \sqcup K_{\frac{n - 1}{2}}),$$
since $\Delta[\{e\}] \cong K_1 \cong \Delta[\{b\}]$ and $\Delta[S \setminus \{e\}])  \cong K_{\frac{n - 1}{2}}$.

If $n$ is even then the conjugacy classes and centralizers of elements of $D_{2n}$ are given by
\begin{align*}
	&e^{D_{2n}} = \{e\}, (a^{\frac{n}{2}})^{D_{2n}} = \{a^{\frac{n}{2}}\},  
	(a^j)^{D_{2n}} = \{a^j, a^{n - j}\} \text{ for } 1 \leq j \leq \frac{n}{2} -1,\\ 
	&b^{D_{2n}} = \{b, a^2b, a^4b, \dots, a^{n - 2}b\}, (ab)^{D_{2n}} = \{ab, a^3b, a^5b, \dots, a^{n - 1}b\}
\end{align*}
and  
\begin{align*}
	C_{D_{2n}}(e) = &D_{2n} = C_{D_{2n}}(a^{\frac{n}{2}}), C_{D_{2n}}(a^i) = \langle a \rangle \text{ for } 1 \leq i \leq n - 1 \text{ and } i \ne \frac{n}{2}\\  
	&C_{D_{2n}}(a^ib) =\{e, a^{\frac{n}{2}},  a^ib, a^{i + \frac{n}{2}}b\} \text{ for } 1 \leq i \leq n, 
\end{align*}
\noindent respectively. It can be seen that 
$$\mathrm{CSCom}(D_{2n}) \cong \Delta[K_1, K_1, \underbrace{K_2, \dots, K_2}_{(\frac{n}{2}-1)\text{-times}}, K_{\frac{n}{2}}, K_{\frac{n}{2}}],$$ 
where $V(\Delta) = \{e, a^{\frac{n}{2}}, a,  \dots, a^{\frac{n}{2}-1}, b, ab\}$ which is identified with the set $\{1, 2, \dots,$  $\frac{n}{2} + 3\}$ preserving the order. Note that $1$ and $a^{\frac{n}{2}}$ are the only dominant vertices of $\Delta$.
The subgraph $\Delta[S]$ is complete if $S=\{e, a^{\frac{n}{2}}, a,  \dots, a^{\frac{n}{2} - 1}\}$.

If $\frac{n}{2}$ is even then $\Delta[\{b, ab\}]$ is isomorphic to $K_1 \sqcup K_1$.
Therefore 
$$\Delta = \Delta[\{e, a^{\frac{n}{2}}\}] \vee (\Delta[\{b, ab\}] \sqcup \Delta[S \setminus \{e, a^{\frac{n}{2}}\}]) \cong K_2 \vee (K_1 \sqcup K_1 \sqcup K_{\frac{n}{2} - 1}),$$
since $\Delta[\{e, a^{\frac{n}{2}}\}] \cong K_2$ and $\Delta[S \setminus \{e, a^{\frac{n}{2}}\}])  \cong K_{\frac{n}{2} -1}$.

If $\frac{n}{2}$ is odd then $\Delta[\{b, ab\}]$ is a complete graph and it is isomorphic to $K_2$.
Therefore 
$$\Delta = \Delta[\{e, a^{\frac{n}{2}}\}] \vee (\Delta[\{b, ab\}] \sqcup \Delta[S \setminus \{e, a^{\frac{n}{2}}\}]) \cong K_2 \vee (K_2 \sqcup K_{\frac{n}{2} - 1}).$$
Hence, the result follows.	
\end{proof}

\begin{theorem}\label{CSCom-Q4n}
	Let  $Q_{4n} = \langle a, b : a^{2n} = e, a^n =b^2, bab^{-1} = a^{-1}\rangle$ be the generalized quaternion group. Then
	\[
	\mathrm{C SCom}(Q_{4n}) \cong   \begin{cases}
		K_2 \vee (K_1 \sqcup K_1 \sqcup K_{n - 1})[K_1, K_1, \underbrace{K_2, \dots, K_2}_{(n-1)\text{-times}}, K_{n}, K_{n}], &	\text{ if $n$ is even}	\\
		K_2 \vee (K_2 \sqcup K_{n - 1})[K_1, K_1, \underbrace{K_2, \dots, K_2}_{(n-1)\text{-times}}, K_{n}, K_{n}], &	\text{ if $n$ is odd.}
		\end{cases}
	\]
\end{theorem}
\begin{proof}
The conjugacy classes and centralizers of elements of $Q_{4n}$ are given by
\begin{align*}
	&e^{Q_{4n}} = \{e\}, (a^n)^{Q_{2n}} = \{a^n\},  
	(a^i)^{Q_{4n}} = \{a^i, a^{2n - i}\} \text{ for } 1 \leq i \leq n-1,\\ 
	&b^{Q_{4n}} = \{b, a^2b, a^4b, \dots, a^{2n - 2}b\}, (ab)^{Q_{4n}} = \{ab, a^3b, a^5b, \dots, a^{2n - 1}b\}
\end{align*}
and  
\begin{align*}
	C_{Q_{4n}}(e) = &Q_{4n} = C_{Q_{4n}}(a^n), C_{Q_{4n}}(a^i) = \langle a \rangle \text{ for } 1 \leq i \leq 2n - 1 \text{ and } i \ne n\\  
	&C_{Q_{4n}}(a^ib) =\{e, a^n,  a^ib, a^{i + n}b\} \text{ for } 1 \leq i \leq 2n, 
\end{align*}
\noindent respectively. It can be seen that 
$$\mathrm{CSCom}(Q_{4n}) \cong \Delta[K_1, K_1, \underbrace{K_2, \dots, K_2}_{(n-1)\text{-times}}, K_{n}, K_{n}],$$ 
where $V(\Delta) = \{e, a^n, a,  \dots, a^{n-1}, b, ab\}$ which is identified with the set $\{1, 2,$ $ \dots,  n + 3\}$ preserving the order. Note that $e$ and $a^n$ are the only dominant vertices of $\Delta$.
The subgraph $\Delta[S]$ is complete if $S=\{e, a^n, a,  \dots, a^{n - 1}\}$.

If $n$ is even then $\Delta[\{b, ab\}]$ is isomorphic to $K_1 \sqcup K_1$.
Therefore 
$$\Delta = \Delta[\{e, a^n\}] \vee (\Delta[\{b, ab\}] \sqcup \Delta[S \setminus \{e, a^n\}]) \cong K_2 \vee (K_1 \sqcup K_1 \sqcup K_{n - 1}),$$
since $\Delta[\{e, a^n\}] \cong K_2$ and $\Delta[S \setminus \{e, a^n\}])  \cong K_{n -1}$.

If $n$ is odd then $\Delta[\{b, ab\}]$ is a complete graph and it is isomorphic to $K_2$.
Therefore 
$$\Delta = \Delta[\{e, a^n\}] \vee (\Delta[\{b, ab\}] \sqcup \Delta[S \setminus \{e, a^n\}]) \cong K_2 \vee (K_2 \sqcup K_{n - 1}).$$	
Hence, the result follows.
\end{proof}
\section{Wiener index of supergraphs}

In this section we use the representation of the supergraphs
as generalized compositions to derive some of their properties. 


\paragraph{Definition} \cite{w,survey, survey1}
Let $\Gamma$ be a graph and let $u$ and $v$ be two vertices of $\Gamma$. The distance between $u$ and $v$, denoted by $d_\Gamma(u,v)$, is defined to be the length of the shortest path between $u$ and $v$. The \emph{ Wiener index} of the graph $\Gamma$, denoted by $W(\Gamma)$,  is defined to be the sum of all distances between any two vertices of $\Gamma$.
\medskip

Let $\mathcal{G}$ be a graph such that $\mathcal{G} = \mathcal{H}[\Gamma_1,\dots,\Gamma_k]$ where $\mathcal{H}$ is a connected graph on $k$ vertices and the factor graphs $\Gamma_i$ are either complete or empty (graph with vertices but no edges). Let $V(\mathcal{H}) =\{1,2,\ldots,k\}$ and $V(\Gamma_i) = \{v_i^{1},\dots,v_i^{n_i}\}$ be the vertex set of the graph $\Gamma_i$, then $V(\mathcal{G}) = \sqcup_{i=1}^kV(\Gamma_i)$. 
The following lemma from \cite{sga} explains the distance between any two vertices of $\mathcal{G}$.

\begin{lemma}\label{general distance}
Let $v_i^p$ and $v_j^q$ ($1 \le p \le n_i \mbox{ and }1 \le q \le n_j$) be two arbitrary vertices of $\mathcal{G}$. Then
\begin{equation*}
d_\mathcal{G}(v_i^p,v_j^q) = \begin{cases}
		1, &\mbox{ if } i=j \mbox{ and }\Gamma_i \mbox{ is complete}\\
		2, &\mbox{ if } i=j \mbox{ and }\Gamma_i \mbox{ is empty}\\
		d_\mathcal{H}(i,j), &\mbox{ if } i \ne j. 
	\end{cases}
\end{equation*}
\end{lemma}

\begin{proof}
	Suppose $i=j$ and $\Gamma_i$ is an empty graph. Then $v_i^p$ and $v_i^q$ are not adjacent in $\Gamma_i$. Now, since the graph $\mathcal{H}$ is connected, there exists a vertex $t$ in $\mathcal{H}$ which is adjacent to $i$ and every vertex of $\Gamma_i$ is adjacent to every vertex of $\Gamma_t$ in the graph $\mathcal{G}$. Therefore, if $s$ is an vertex of $\Gamma_t$ then $v_i^p \rightarrow s \rightarrow v_i^q$ is a path of length two in $\mathcal{G}$. 
	
	Suppose $i \ne j$. Then $v_i^p$ and $v_j^q$ are in two distinct factors $\Gamma_i$ and $\Gamma_j$ of $\mathcal{G}$. Since, the graph $\mathcal{H}$ is connected, there exists a path $i=i_0 \rightarrow i_1 \rightarrow \cdots \rightarrow i_k = j$ of length $d_\mathcal{H}(i,j)$ in $\mathcal{H}$. Now,   $v_{i}^p \rightarrow v_{i_1}^{p_1} \rightarrow \cdots \rightarrow v_{i_k}^{q}$ where $1 \le p_t \le n_{i_t}$ is a path of length $d_\mathcal{H}(i,j)$ joining $v_i^p$ and $v_j^q$. From the definition of graph $\mathcal{G}$, we see that this path has to be the minimal one between $v_i^p$ and $v_i^q$. This completes the proof of the lemma. 
\end{proof}

From the above lemma, we get the following result.

\begin{lemma}\label{mainlem}
	Let $\mathcal{G}$ be the graph defined above. Then the  Wiener index of $\mathcal{G}$ is given by \begin{equation*}
		W(G) = \sum_{\substack{1 \le i \le k \\ \Gamma_i \hbox{ is complete}}}\binom{n_i}{2} + \sum_{\substack{1 \le i \le k \\ \Gamma_i \hbox{ is empty}}}2\binom{n_i}{2}+\sum_{1 \le i<j \le k}n_in_jd_\mathcal{H}(i,j),
	\end{equation*} where $|V(\Gamma_i)| = n_i$. 
\end{lemma}



Let $A$ be a finite simple graph with vertex set $G$. Let $B = \{V_1,\dots,V_k\}$ be a partition of the vertex set $G$ and $\sim$ is the corresponding equivalence relation. Let $\Gamma_A^B$ be the associated $\partype$ super$\grtype$ graph. Let $R = \{v_1,\dots,v_k\}$ be a set of representatives of the classes $ \{V_1,\dots,V_k\}$. 
%
Then, by Proposition \ref{p:join},
$\Gamma_A^B = \Delta[K_{n_1}, \dots, K_{n_k}]$, where  $n_i = |V_i|\geq 1$ for $1 \leq i \leq k$.
 Hence, we can use Lemma \ref{mainlem} to calculate the Wiener index of $\Gamma_A^B$. 
\begin{theorem}\label{W-Index}
The Wiener index of the graph  $\Gamma_A^B = \Delta[K_{n_1}, \dots, K_{n_k}]$, where $V(\Delta) = \{v_1, \dots, v_k\}$, is equal to
\[
W(\Gamma_A^B) = \sum_{\substack{1 \le i \le k \\ n_i > 1}}\binom{n_i}{2} + \sum_{1 \le i<j \le k}n_in_jd_{\Delta}(v_i,v_j).
\]
\end{theorem}

\begin{proof}
        The result follows from Lemma \ref{mainlem}   noting that  $\underset{{\substack{1 \le i \le k \\ n_i = 1}}}{\sum}\binom{n_i}{2} = 0$.
\end{proof}

In the following results, we will calculate $W(\Gamma_A^B)$  when $\Gamma_A^B$ is equality superA graph and conjugacy superA graph  for dihedral groups and quaternion group where the graph A is  commuting graph.

\begin{theorem}\label{WI-CSCom-D2n}
	Let   $D_{2m} = \langle a, b : a^m =b^2= e, bab^{-1} = a^{-1}\rangle$ the dihedral group. Then
	\[
	W(\mathrm{CSCom}(D_{2m})) = \begin{cases}
		\frac{1}{2}(5m^2  - 2m + 3), &\text{ if  $m$ is odd}\\
		2m^2  + 7m -20, &	\text{ if $m$ and $\frac{m}{2}$ are even}\\
		2m^2  + 7m -24, &	\text{ if $m$ is even and $\frac{m}{2}$ is odd.}
	\end{cases}
	\]
\end{theorem}
\begin{proof}
	If $m$ is odd then, by Theorem 	\ref{CSCom-D2n}, we have 
	$$
	\mathrm{CSCom}(D_{2m}) = {\mathcal{G}}_1 \vee ({\mathcal{G}}_2 \sqcup {\mathcal{G}}_3)[{\mathcal{G}}_4, \underbrace{{\mathcal{G}}_5, \dots, {\mathcal{G}}_5}_{(\frac{m-1}{2})\text{-times}}, {\mathcal{G}}_6],
	$$ 
	where ${\mathcal{G}}_1 \cong {\mathcal{G}}_2 \cong {\mathcal{G}}_4 \cong K_1$, ${\mathcal{G}}_3 \cong K_{\frac{m - 1}{2}}$, ${\mathcal{G}}_5 \cong K_2$ and ${\mathcal{G}}_6 \cong K_m$. Let $V({\mathcal{G}}_1) = \{v_1\}$, $V({\mathcal{G}}_2) = \{v_2\}$, $V({\mathcal{G}}_3) = \{v_3, \dots, v_{\frac{m-1}{2} + 2}\}$ and $\Delta = {\mathcal{G}}_1 \vee ({\mathcal{G}}_2 \sqcup {\mathcal{G}}_3)$. In view of Theorem \ref{W-Index}, we have $n_1 = 1, n_2 = \cdots = n_{\frac{m-1}{2} + 1} = 2$ and $n_{\frac{m-1}{2} + 2} = m$. It follows that 
\[
d_{\Delta}(v_1, v_j) = 1 \text{ for } 2 \leq j \leq \frac{m-1}{2} + 2, \, d_{\Delta}(v_2, v_j) = 2 \text{ for } 3 \leq j \leq \frac{m-1}{2} + 2
\]
\[ 
\text{ and } d_{\Delta}(v_i, v_j) = 1 \text{ for } 3 \leq i < j \leq \frac{m-1}{2} + 2.
\]
Also,
\[
n_1n_j = 2 \text{ for } 2 \leq j \leq \frac{m-1}{2} + 1, n_1n_{\frac{m-1}{2} + 2} = m, 
n_2n_j = 4 \text{ for } 3 \leq j \leq \frac{m-1}{2} + 1, n_2n_{\frac{m-1}{2} + 2} = 2m 
\]
and
\[
n_in_j = 4 \text{ for } 3 \leq i < j \leq \frac{m-1}{2} + 1, n_{\frac{m-1}{2} + 1}n_{\frac{m-1}{2} + 2} = 2m.
\]
Therefore,
\begin{align*}
\sum_{1 \le i<j \le \frac{m-1}{2} + 2}n_in_jd_{\Delta}(v_i,v_j) =& \sum_{2 \le j \le \frac{m-1}{2} + 2}n_1n_jd_{\Delta}(v_1,v_j) + \sum_{3 \le j \le \frac{m-1}{2} + 2}n_2n_jd_{\Delta}(v_2,v_j)\\
& \qquad\qquad\qquad + \sum_{3 \le i<j \le \frac{m-1}{2} + 2}n_in_jd_{\Delta}(v_i,v_j)\\
=& 2m - 1 + 8m -12 + (m-2)(2m-5)\\
=& 2m^2 -m + 2
\end{align*}
and 
\begin{align*}
\sum_{\substack{1 \le i \le \frac{m-1}{2} + 2 \\ n_i > 1}}\binom{n_i}{2} = \underbrace{1 + \cdots + 1}_{\text{$(\frac{m-1}{2})$-times}} + \frac{m(m-1)}{2} = \frac{m^2-1}{2}.
\end{align*}	
Hence, by Theorem \ref{W-Index}, we have
\[
W(\mathrm{CSCom}(D_{2m})) = \frac{m^2-1}{2} + 2m^2 -m + 2 = \frac{1}{2}(5m^2  - 2m + 3).
\]

If $m$ and $\frac{m}{2}$ both are even then, by Theorem 	\ref{CSCom-D2n}, we have 
$$
\mathrm{CSCom}(D_{2m}) = {\mathcal{G}}_1 \vee ({\mathcal{G}}_2 \sqcup {\mathcal{G}}_3 \sqcup {\mathcal{G}}_4)[{\mathcal{G}}_5, {\mathcal{G}}_6, \underbrace{{\mathcal{G}}_7, \dots, {\mathcal{G}}_7}_{(\frac{m}{2} - 1)\text{-times}}, {\mathcal{G}}_8, {\mathcal{G}}_9],
$$ 
where ${\mathcal{G}}_1 \cong {\mathcal{G}}_7 \cong K_2$, ${\mathcal{G}}_2 \cong {\mathcal{G}}_3 \cong {\mathcal{G}}_5 \cong {\mathcal{G}}_6 \cong K_1$, ${\mathcal{G}}_4 \cong K_{\frac{m}{2} - 1}$,  and ${\mathcal{G}}_8 \cong {\mathcal{G}}_9 \cong K_{\frac{m}{2}}$. Let $V({\mathcal{G}}_1) = \{v_1, v_2\}$, $V({\mathcal{G}}_2) = \{v_3\}$, $V({\mathcal{G}}_3) = \{v_4\}$, $V({\mathcal{G}}_4) = \{v_5, \dots, v_{\frac{m}{2} + 3}\}$ and $\Delta = {\mathcal{G}}_1 \vee ({\mathcal{G}}_2 \sqcup {\mathcal{G}}_3) \sqcup {\mathcal{G}}_4$. In view of Theorem \ref{W-Index}, we have $n_1 = n_2 = 1, n_3 = \cdots = n_{\frac{m}{2} + 1} = 2$ and $n_{\frac{m}{2} + 2} = n_{\frac{m}{2} + 3} = \frac{m}{2}$. It follows that 
\[
d_{\Delta}(v_1, v_j) = 1 \text{ for } 2 \leq j \leq \frac{m}{2} + 3, \, d_{\Delta}(v_2, v_j) = 1 \text{ for } 3 \leq j \leq \frac{m}{2} + 3,
\]
\[
d_{\Delta}(v_3, v_j) = 2 \text{ for } 4 \leq j \leq \frac{m}{2} + 3, \, d_{\Delta}(v_4, v_j) = 2 \text{ for } 5 \leq j \leq \frac{m}{2} + 3
\]
\[ 
\text{ and } d_{\Delta}(v_i, v_j) = 1 \text{ for } 5 \leq i < j \leq \frac{m}{2} + 3.
\]
Also,
\[
n_1n_2 = 1, n_1n_j = 2 \text{ for } 3 \leq j \leq \frac{m}{2} + 1, \, n_1n_{\frac{m}{2} + 2} = n_1n_{\frac{m}{2} + 3}  = \frac{m}{2},
\]
\[
n_2n_j = 2 \text{ for } 3 \leq j \leq \frac{m}{2} + 1, \, n_2n_{\frac{m}{2} + 2} = n_2n_{\frac{m}{2} + 3}  = \frac{m}{2}, 
\]
\[
n_in_j = 4 \text{ for } 3 \leq i < j \leq \frac{m}{2} + 1, \, 
n_in_{\frac{m}{2} + 2} = n_in_{\frac{m}{2} + 3}  = m \text{ for } 3 \leq i  \leq \frac{m}{2} + 1
\]
and
\[
n_{\frac{m}{2} + 2}n_{\frac{m}{2} + 3} = \frac{m^2}{4}.
\]
Therefore,
\begin{align*}
	\sum_{1 \le i<j \le \frac{m}{2} + 3}n_in_jd_{\Delta}(v_i,v_j) =& \sum_{2 \le j \le \frac{m}{2} + 3}n_1n_jd_{\Delta}(v_1,v_j) + \sum_{3 \le j \le \frac{m}{2} + 3}n_2n_jd_{\Delta}(v_2,v_j)\\
	& + \sum_{4 \le j \le \frac{m}{2} + 3}n_3n_jd_{\Delta}(v_3,v_j) 
	+ \sum_{5 \le j \le \frac{m}{2} + 3}n_4n_jd_{\Delta}(v_4,v_j)\\
	& \qquad\qquad\qquad + \sum_{5 \le i<j \le \frac{m}{2} + 3}n_in_jd_{\Delta}(v_i,v_j)\\
	=& (2m - 1) +  (2m - 2)	+ (8m - 16) + (8m - 24) + (\frac{7m^2}{4} - 13m + 24)\\
	=&\frac{7m^2}{4} + 7m - 19
\end{align*}
and 
\begin{align*}
	\sum_{\substack{1 \le i \le \frac{m}{2} + 3 \\ n_i > 1}}\binom{n_i}{2} = \underbrace{1 + \cdots + 1}_{\text{$(\frac{m}{2} - 1)$-times}} \,\, + \,\, 2\times \frac{m(m-2)}{8} = \frac{m^2}{4} - 1.
\end{align*}	
Hence, by Theorem \ref{W-Index}, we have
\[
W(\mathrm{CSCom}(D_{2m})) = \frac{m^2}{4} - 1 + \frac{7m^2}{4} + 7m - 19 = 2m^2  + 7m -20.
\]

If $m$ is even and $\frac{m}{2}$ is odd then, by Theorem 	\ref{CSCom-D2n}, we have 
$$
\mathrm{CSCom}(D_{2m}) = {\mathcal{G}}_1 \vee ({\mathcal{G}}_2  \sqcup {\mathcal{G}}_3)[{\mathcal{G}}_4, {\mathcal{G}}_5, \underbrace{{\mathcal{G}}_6, \dots, {\mathcal{G}}_6}_{(\frac{m}{2} - 1)\text{-times}}, {\mathcal{G}}_7, {\mathcal{G}}_8],
$$ 
where ${\mathcal{G}}_1 \cong {\mathcal{G}}_2 = {\mathcal{G}}_6 \cong K_2$, ${\mathcal{G}}_3 \cong K_{\frac{m}{2} - 1}$, ${\mathcal{G}}_4 \cong {\mathcal{G}}_5 \cong K_1$   and ${\mathcal{G}}_7 \cong {\mathcal{G}}_8 \cong K_{\frac{m}{2}}$. Let $V({\mathcal{G}}_1) = \{v_1, v_2\}$, $V({\mathcal{G}}_2) = \{v_3, v_4\}$, $V({\mathcal{G}}_3) = \{v_5, \dots, v_{\frac{m}{2} + 3}\}$ and $\Delta = {\mathcal{G}}_1 \vee ({\mathcal{G}}_2 \sqcup {\mathcal{G}}_3)$. In view of Theorem \ref{W-Index}, we have $n_1 = n_2 = 1, n_3 = \cdots = n_{\frac{m}{2} + 1} = 2$ and $n_{\frac{m}{2} + 2} = n_{\frac{m}{2} + 3} = \frac{m}{2}$. It is easy to see that $d_{\Delta}(v_3,v_4) = 1$ and all other values of $d_{\Delta}(v_i,v_j)$, $n_in_j$ and $\sum_{\substack{1 \le i \le \frac{m}{2} + 3 \\ n_i > 1}}\binom{n_i}{2}$  are same as the values described in the case when $m$ and $\frac{m}{2}$ both are even. Therefore, 
$\sum_{4 \le j \le \frac{m}{2} + 3}n_3n_jd_{\Delta}(v_3,v_j) = 4 + 8\times (\frac{m}{2} - 3) + 4m =  8m - 20$ and so
\begin{align*}
\sum_{1 \le i<j \le \frac{m}{2} + 3}n_in_jd_{\Delta}(v_i,v_j) &= (2m - 1) +  (2m - 2)	+ (8m - 20) + (8m - 24) + (\frac{7m^2}{4} - 13m + 24)\\
& = \frac{7m^2}{4} + 7m - 23.
\end{align*}
Hence, by Theorem \ref{W-Index}, we have
\[
W(\mathrm{CSCom}(D_{2m})) = \frac{m^2}{4} - 1 + \frac{7m^2}{4} + 7m - 23 = 2m^2  + 7m -24.
\]
This completes the proof.
\end{proof}

\begin{theorem}\label{WI-CSCom-Q4n}
	Let  $Q_{4n} = \langle a, b : a^{2n} = e, a^n =b^2, bab^{-1} = a^{-1}\rangle$ the generalized quaternion group. Then
	\[
	W(\mathrm{CSCom}(Q_{4n})) = \begin{cases}
		8n^2 + 14n - 20, &	\text{ if $n$ is even}	\\
		8n^2 + 14n - 24, &	\text{ if $n$ is odd.}
	\end{cases}
	\]	
\end{theorem}
\begin{proof}
By Theorems \ref{CSCom-D2n} and \ref{CSCom-Q4n} we have 	$\mathrm{CSCom}(Q_{4n}) \cong \mathrm{CSCom}(D_{2\times 2n})$. Therefore, if $n$ is even then using Theorem \ref{WI-CSCom-D2n} we get
\[
W(\mathrm{CSCom}(Q_{4n})) = \begin{cases}
	2(2n)^2  + 7\times 2n -20, & \text{ if $n$ is even}\\
    2(2n)^2  + 7\times 2n -24, & \text{ if $n$ is odd}.
\end{cases}
\]
Hence, the result follows.
\end{proof}

\begin{theorem}\label{WI-ESCom-D2n}
	Let   $D_{2n} = \langle a, b : a^n =b^2= e, bab^{-1} = a^{-1}\rangle$ the dihedral group. Then
	\[
	W(\mathrm{ESCom}(D_{2n})) = \begin{cases}
	\frac{1}{2}(7n^2 - 5n), &\text{ if  $n$ is odd}\\
	\frac{1}{2}(7n^2 - 8n), &	\text{ if $n$ is even.}
	\end{cases}
	\]
\end{theorem}
\begin{proof}
If $n$ is odd then, by Theorem 	\ref{ESCom-D2n}, we have 
$$
\mathrm{ESCom}(D_{2n}) = {\mathcal{G}}_1 \vee ({\mathcal{G}}_2 \sqcup \cdots \sqcup {\mathcal{G}}_{n+1} \sqcup {\mathcal{G}}_{n+2}) :=\mathcal{G},
$$ 
where ${\mathcal{G}}_i \cong   K_1$ for $1\leq i \leq n+1$ and ${\mathcal{G}}_{n+2} \cong K_{n - 1}$. Let $V({\mathcal{G}}_i) = \{v_i\}$ for $1\leq i \leq n+1$ and  $V({\mathcal{G}}_{n+2}) = \{v_{n+2}, \dots, v_{2n}\}$. Therefore,
\[
d_{\mathcal{G}}(v_1, v_j) = 1 \text{ for } 2 \leq j \leq 2n, d_{\mathcal{G}}(v_i, v_j) = 2 \text{ for } 2 \leq i < j \leq n + 1
\]
\[ 
\text{ and } d_{\mathcal{G}}(v_i, v_j) = 1 \text{ for } n+2 \leq i < j \leq 2n.
\]
Therefore,
\begin{align*}
		W(\mathrm{ESCom}(D_{2n})) & = \sum_{1 \le i<j \le 2n}d_{\mathcal{G}}(v_i,v_j)\\
		& = \sum_{2 \le j \le 2n}d_{\mathcal{G}}(v_1,v_j) + \sum_{2 \le i < j \le n + 2}d_{\mathcal{G}}(v_i,v_j) + \sum_{n+2 \le i< j \le 2n}d_{\mathcal{G}}(v_i,v_j)\\
		& = (2n - 1) + 3n(n - 1) + \frac{n^2 - 3n + 2}{2}\\
		& = \frac{n(7n - 5)}{2}.
\end{align*}

If $n$ is even then, by Theorem 	\ref{ESCom-D2n}, we have 
$$
\mathrm{ESCom}(D_{2n}) = {\mathcal{G}}_1 \vee ({\mathcal{G}}_2 \sqcup \cdots \sqcup {\mathcal{G}}_{\frac{n}{2}+1} \sqcup {\mathcal{G}}_{\frac{n}{2}+2}) := \mathcal{G},
$$ 
where ${\mathcal{G}}_i \cong   K_2$ for $1\leq i \leq \frac{n}{2}+1$ and ${\mathcal{G}}_{\frac{n}{2}+2} \cong K_{n - 2}$. Let $V({\mathcal{G}}_i) = \{v_i, v_{i + 1}\}$ for $1\leq i \leq \frac{n}{2}+1$ and  $V({\mathcal{G}}_{n+2}) = \{v_{\frac{n}{2}+3}, \dots, v_{2n}\}$. Therefore,
\[
d_{\mathcal{G}}(v_1, v_j) = 1 \text{ for } 2 \leq j \leq 2n, d_{\mathcal{G}}(v_2, v_j) = 1 \text{ for } 3 \leq j \leq 2n,
\]
\[
 d_{\mathcal{G}}(v_i, v_j) = 2 \text{ for } 3 \leq i < j \leq \frac{n}{2} + 2
\text{ and } d_{\mathcal{G}}(v_i, v_j) = 1 \text{ for } \frac{n}{2}+3 \leq i < j \leq 2n.
\]
Therefore,
\begin{align*}
	W(\mathrm{ESCom}&(D_{2n}))  = \sum_{1 \le i<j \le 2n}d_{\mathcal{G}}(v_i,v_j)\\
	& = \sum_{2 \le j \le 2n}d_{\mathcal{G}}(v_1,v_j) + \sum_{3 \le j \le 2n}d_{\mathcal{G}}(v_2,v_j)+ \sum_{3 \le i < j \le \frac{n}{2} + 2}d_{\mathcal{G}}(v_i,v_j) + \sum_{\frac{n}{2}+2 \le i< j \le 2n}d_{\mathcal{G}}(v_i,v_j)\\
	& = (2n - 1) + (2n - 2) + \frac{6n^2 -11n}{2} + \frac{n^2 - 5n + 6}{2}\\
	& = \frac{n(7n - 8)}{2}.
\end{align*}	
Hence, the result follows.
\end{proof}

\begin{theorem}\label{WI-ESCom-Q4n}
	Let  $Q_{4n} = \langle a, b : a^{2n} = e, a^n =b^2, bab^{-1} = a^{-1}\rangle$ the generalized quaternion group. Then
	\[
	W(\mathrm{ESCom}(Q_{4n})) = 14n^2 - 8n.
	\]
\end{theorem}
\begin{proof}
By Theorems \ref{ESCom-D2n} and \ref{ESCom-Q4n} it follows that 
$\mathrm{ESCom}(Q_{4n}) \cong \mathrm{ESCom}(D_{2\times 2n})$. Therefore, by Theorem \ref{WI-ESCom-D2n}, we get $W(\mathrm{ESCom}(Q_{4n})) = \frac{1}{2}(7\times 4n^2 - 8 \times 2n) = 14n^2 - 8n$.
\end{proof}

\end{document}